\documentclass[a4paper,11pt]{article}
\usepackage{amsfonts,amssymb,latexsym,amsmath}
\usepackage{multicol}
\newcounter{num}[section]
\newcommand{\Num}{\refstepcounter{num}%
\textbf{\arabic{section}.\arabic{num}}}

\newcommand{\Theorem}{\textbf{Theorem~}}
\newcommand{\Proof}{\textbf{Proof}}
\newcommand{\Conj}{\textbf{Conjecture~}}
\newcommand{\Lemma}{\textbf{Lemma~}}

\newcommand{\Prop}{\textbf{Proposition~}}

\newcommand{\De}{\Delta}
\newcommand{\de}{\delta}
\newcommand{\UT}{{\mathrm{UT}(n,\Rb)}}
\newcommand{\UTk}{{\mathrm{UT}(k,\Rb)}}
\newcommand{\Tr}{{\mathrm{Tr}}}
\newcommand{\dpi}{{2\pi i}}
\newcommand{\ct}{\tilde{c}}
\newcommand{\ga}{\gamma}

\newcommand{\la}{\lambda}
\newcommand{\La}{\Lambda}
\newcommand{\vphi}{{\varphi}}
\newcommand{\hvphi}{{\hat{\vphi}}}
\newcommand{\ut}{{\ux\tx(n,\Rb)}}
\newcommand{\Xx}{{\cal X}}

\newcommand{\gx}{{\mathfrak g}}
\newcommand{\ux}{{\mathfrak u}}
\newcommand{\tx}{{\mathfrak t}}
\newcommand{\px}{{\mathfrak p}}
\newcommand{\Ad}{{\mathrm{Ad}}}
\newcommand{\Ann}{\mathrm{Ann}}
\newcommand{\diag}{\mathrm{diag}}

\newcommand{\Rb}{{\Bbb R}}

\renewcommand{\leq}{\leqslant}

\begin{document}
\Large

\title{Subregular characters of the group $\UT$}
\author{A.N.Panov\thanks{Work of the first author is  supported by the RFBR-grants 12-01-00070,
12-01-00137, 13-01-97000-Volga region-a} \and E.V.Surai}

\date{}

\maketitle

\section{\textsc{Introduction}}
Notion of character of  representation plays an important role in
the representation theory of groups. For  finite dimensional
representation a character $\chi(g)$ is defined as a trace of matrix
of  operator  $T_g$. This definition cannot be directly applied for
infinite dimensional  representations. For some representations of
Lie groups one can define a character as a generalized function on
the group as follows.

Extend a representation   $T_g$ of a group $G$ to the representation
of the group algebra  $L^1(G)$ by the formula
$$T_\vphi =  \int\vphi(g)T_g dg, $$
where $dg$ is the left-invariant measure on the group. Suppose that
for any finite function the operator $T_\vphi$ has the trace. Then
the formula
$$(\chi,\vphi) = \Tr(T_\vphi)$$ defines the trace  $\chi(g)$ as a generalized function on the Lie group  $G$.

It was proved in the paper  \cite{K1} that any irreducible
representation of a connected nilpotent Lie group has a trace in
mentioned sense. For any connected nilpotent Lie group and  for a
character of  irreducible representation, which is associated with a
coadjoint orbit  $\Omega$, the formula of A.A.Kirillov is valid
\begin{equation}\label{int}
 (\chi,\vphi) = \int \hvphi(a) d_\Omega\mu(a), \end{equation}
 where $$ \hvphi(a)
= \int \vphi(\exp(x))e^{2\pi i (a,x)} dx$$ is the  Fourier
transformation with respect to a given Lebesgue measure  $dx$ on the
Lie algebra  $\gx$ of  $G$, ~$d_\Omega\mu$ is  the invariant measure
on the orbit $\Omega$ (see \cite{K1} and \cite[chapter 3]{K2}). But
it may be very difficult to use this formula for calculation of the
character, when the orbit is a complicated manifold, defined by a
system of many equations. It will be better to use the other
approach:  if we represent the operator $T_\vphi$ in the integral
form
\begin{equation}\label{iint}T_\vphi
F(x)=\int A(x,y)F(y)dy,\end{equation}
 then the trace is equal to  $
\int A(x,x) dx$ (see \cite{K1}). The $\de$-function type multipliers
rise in the character; it implies that  the support of  character
doesn't coincide with the group. Calculations of characters   for
different connected nilpotent Lie groups imply the following general
conjecture.

\Conj~\Num\label{hy}. Given a connected nilpotent Lie group,  the
support of a character of irreducible representation, associated
with an coadjoint orbit  $\Omega$,  coincides with  the closure of
union of stabilizers  $G^f$, where $f\in \Omega$.

  The regular characters  (i.e. characters of irreducible representations, associated with
  orbits of maximal   dimension) of $\UT$    were calculated in the paper  ~\cite{K1}.
Recall that any  coadjoint orbit has even dimension. An coadjoint
orbit is called  subregular if its dimension equals to $d-2$, where
$d$ is a maximal dimension of coadjoint orbits. A character of
irreducible representation, associated with a subregular coadjoint
orbit, is called a subregular character. Subregular coadjoint orbits
of the group $\UT$ were classified in ~\cite{IP}. Basing on this
classification we obtain formulas (theorems \ref{subeven} and
\ref{subodd}) for subregular characters. In this paper we also
calculate regular characters  (see theorems  \ref{theven} and
\ref{thodd}), since  there is no proof of this formulas in
~\cite{K1} and there is some unexpected misprint (in the case of odd
$n$). The obtained formulas confirm Conjecture  \ref{hy} in the case
of regular and subregular characters of the group $\UT$ (see Theorem
 \ref{thhy}).

Note that subregular characters for the case of  finite field were
obtained in  \cite{Ig}. We use the following notations:\\
1) ~ if $C$ is a matrix with functional entries, then $\de(C)$ is
a product of $\de$-functions at zero of its entries ;\\
2)~ if $C$ is an unitriangular matrix with functional entries upper
the diagonal, then we preserve the notation
 $\de(C)$ for the product of $\de$-functions at zero of  its entries upper the diagonal;\\
3) ~ if $C_1,\ldots, C_m$ is a system of matrix, then
$\de(C_1,\ldots, C_m) = \de(C_1)\cdots\de(C_m)$.

\section{\textsc{Regular characters}}

The group  $G= \UT$ is a group of upper triangular matrices of size
 $n$ with units on the diagonal. Its Lie algebra  $\gx = \ut$ consists of all
upper triangular matrices with zeros on the diagonal. Applying the
Killing form  $(\cdot,\cdot)$ we identify  the conjugate space
$\gx^*$ with the set of all lower triangular matrices with seros on
the diagonal. According to the orbit method of A.A.Kirillov, there
exists one to one correspondence between irreducible representations
of a connected nilpotent Lie group and its coadjoint orbits. Any
element  $f\in\gx^*$ has a polarization  $\px$ (i.e. a subalgebra
that is a maximal izotropic subspace with respect to the skew
symmetric bilinear form  $f([x,y])$). The irreducible
representation, associated with the coadjoint orbit  $\Omega(f)$, is
induced from one dimensional representation
$$\xi(\exp(x))
= e^{\dpi f(x)},~ \mbox{where}\quad x\in\px, $$ of the subgroup  $\exp(\px)$.\\

\subsection{ Case of even  $n$.} Let $n=2k$.
Let us represent the general element  $g\in \UT$ in the form  of
block matrix
\begin{equation}\label{gg} g =
\left(\begin{array}{cc} C_{11}&C_{12}\\
0&C_{22}\end{array}\right),\end{equation} where $C_{ij}$ are blocks
of size  $k\times k$. Any regular coadjoint orbit of the group
$\UT$, where $n=2k$, contains a unique element of the form
\begin{equation}\label{ff}f =
\begin{array}{lll}
\left(\begin{array}{cc} 0&0\\
\La&0\end{array}\right),& \quad\quad \La = \left(\begin{array}{ccc}
0&\ldots&\la_k\\
\vdots&\ddots&\vdots\\
\la_1&\ldots&0\end{array}\right),
\end{array}\end{equation}
where $\la_1,\ldots,\la_{k}\in \Rb$ and $\la_1, \ldots,\la_{k-1}\ne
0$ (see ~\cite[\S 3]{IP}). We use the following notations.  For
every  $1\leq s\leq k$ we denote \\
1)~ $
 \De_s = \left|\begin{array}{ccc} c_{s,1}&\ldots& c_{s,k-s+1}\\
\vdots&\ddots&\vdots\\
c_{k,1}&\ldots & c_{k,k-s+1}\end{array}\right| $ the left lower
minor of the block  $C_{12} = \left(c_{ij}\right)_{i,j=1}^k$,\\
2) ~ $P_s = \frac{(-1)^{k-s}\De_{s}}{\De_{s+1}}$.\\

\Theorem\Num\label{theven} (см. \cite{K1}). Let $n=2k$. The
character of the irreducible representation, associated with the
orbit of element  $f\in \ut^*$ of  (\ref{ff}), has the form
  \begin{equation}\label{regeven} \chi(g) = \de(C_{11}, C_{22})
 \chi_\La^*(C_{12}),\end{equation}

where
 $$\begin{array}{c}\chi_\La^*(C_{12}) = \frac{e^{\dpi(\la_1 P_1+\ldots+\la_k
P_k)}}{\mu_0 \cdot |\De_2 \De_3 \cdots \De_k|}, \quad\quad\quad
 \mu_0 = |\la_1^{k-1}\la_2^{k-2}\cdots \la_k^{0}|.\end{array}$$

\Proof. The subalgebra
$ \left(\begin{array}{cc} 0&*\\
0&0\end{array}\right)$ is a polarization of  $f$. The irreducible
representation, associated with the orbit  $\Omega(f)$, is induced
from the one-dimensional representation  $e^{2\pi i(\La,B)}$ of
subroup
$\left\{ \left(\begin{array}{cc} 0&B\\
0&0\end{array}\right)\right\}$. This representation is realized in
the space  $L^2(\Xx)$, where $\Xx$  consists of matrices  $X=
\diag\left(X_{11}, X_{22}\right)$ with blocks of  $\UTk$, by the
formula
$$
T_gF(X) = e^{\dpi\left(\La,
X_{11}C_{12}(X_{22}C_{22})^{-1}\right)}F(X_{11}C_{11},
X_{22}C_{22}).$$

Operator  $T_g$ is extended to operator
$$T_\vphi F(X) = \int \vphi(C_{11}, C_{12}, C_{22})e^{\dpi\left(\La,
X_{11}C_{12}(X_{22}C_{22})^{-1}\right)}F(X_{11}C_{11},
X_{22}C_{22})dC_{11} dC_{12} dC_{22}.$$ After substitution
$Y_{11}=X_{11}C_{11}$, ~$Y_{22}=X_{22}C_{22}$ we obtain  $T_\vphi$
in the form  (\ref{iint}):
$$ T_\vphi F(X) = \int \vphi(X_{11}^{-1}Y_{11}, C_{12}, X_{22}^{-1}Y_{22})
e^{\dpi\left(\La, X_{11}C_{12}Y_{22}^{-1}\right)}F(Y_{11},
Y_{22})dY_{11} dY_{22} dC_{12}.$$

Then  $$(\chi, \vphi) = \int \vphi(E, C_{12}, E)e^{\dpi\left(\La,
X_{11}C_{12}X_{22}^{-1}\right)}dX_{11} dX_{22} dC_{12}.$$ After
substitution  $X_{22}^{-1}$ by $X_{22}$ we have

\begin{equation}\label{fi}
\chi(g) = \de(C_{11}, C_{22}) \int e^{\dpi\left(\La,
X_{11}C_{12}X_{22}\right)}dX_{11} dX_{22}
\end{equation}

Let    $X_{11}, X_{22}$ be unitriangular matrices of size  $k$ with
entries  $x_{ij}$,~ $y_{ij}$,~ $1\leq i<j\leq k$ respectively.

For every $1\leq s\leq k$ and  $1\leq j\leq k-s+1$ we apply the
notations
$$\ct_{s,j} = c_{s,j} + x_{s,s+1}c_{s+1,j} + \ldots + x_{s,k}c_{k,j},$$
\begin{equation}\label{ss} \chi_s(g) = \int e^{\dpi \la_s\left(\ct_{s,1} y_{1,k-s+1} +
\ldots + \ct_{s,k-s} y_{k-s,k-s+1} + \ct_{s,k-s+1}\right)}dx_s
dy_s,\end{equation}
 where $dx_s = dx_{s,s+1}\ldots dx_{s,k}$ and $dy_s =
dy_{1,k-s+1}\ldots dy_{k-s,k-s+1}$. Taking into account (\ref{ss}),
the formula (\ref{fi}) is rewritten in the form
\begin{equation}\label{cc}
\chi(g) = \de(C_{11}, C_{22})\prod_{s=1}^k \chi_s(g)\end{equation}
Apply  the well  known equality
$$\int e^{\dpi \la(a,x)}dx = \de(a)\frac{1}{|\la|}$$
in  (\ref{ss}),
 we obtain
\begin{equation}\label{sss} \chi_s(g) = \frac{1}{|\la_s|^{k-s}}\int
\de(\ct_{s,1},\ldots,\ct_{s,k-s}) e^{\dpi \la_s
\ct_{s,k-s+1}}dx_s.\end{equation}

 Let
$(x_{s,s+1}^0,\ldots, x_{s,k}^0)$ be a solution of the system of
linear equations
\begin{equation}\label{ee}\left\{\begin{array}{l}
\ct_{s,1} = c_{s,1} + x_{s,s+1}c_{s+1,1} + \ldots + x_{s,k}c_{k,1} = 0,\\
\ldots\\
\ct_{s,k-s} = c_{s,k-s} + x_{s,s+1}c_{s+1,k-s} + \ldots +
x_{s,k}c_{k,k-s} = 0.\end{array}\right.
\end{equation}

Denote $$P'_s =  c_{s,k-s+1} + x_{s,s+1}^0 c_{s+1,k-s+1} + \ldots +
x_{s,k}^0 c_{k,k-s+1}.$$ Let us introduce a new variable $x_{s,s}$.
The vector $(1, x_{s,s+1}^0,\ldots, x_{s,k}^0)$ is a  solution of
the  system of linear equations

\begin{equation}\label{eee}\left\{\begin{array}{l}
x_{s,s} c_{s,1} + x_{s,s+1}c_{s+1,1} + \ldots + x_{s,k}c_{k,1} = 0,\\
\ldots\\
 x_{s,s} c_{s,k-s} + x_{s,s+1}c_{s+1,k-s} + \ldots +
x_{s,k}c_{k,k-s} = 0,\\
 x_{s,s} c_{s,k-s+1} + x_{s,s+1} c_{s+1,k-s+1} +
\ldots + x_{s,k} c_{k,k-s+1} = P'_s.
\end{array}\right.
\end{equation}

Using the Cramer formulas in  (\ref{eee}), we have
$$1 = \frac{(-1)^{k-s}\De_{s+1}P'_s}{\De_s}.$$
It implies $$ P'_s = \frac{(-1)^{k-s}\De_{s}}{\De_{s+1}} =
P_s,\quad\quad\quad \chi_s(g) = \frac{1}{|\la_s|^{k-s} |\De_{s+1}|}
e^{\dpi \la_sP_s}.$$ Substituting  $\chi_s(g)$ into (\ref{cc}), we
obtain the formula for the regular character for even  $n$. $\Box$

\subsection{ Case of odd  $n$.} Let $n=2k+1$.
Let us represent the general element  $g\in \UT$ in the form of
block matrix
\begin{equation}\label{ggg} g =
\left(\begin{array}{ccc} C_{11}&C_{12}&C_{13}\\
0&1&C_{23}\\
0&0&C_{33}
\end{array}\right),\end{equation} where the partition of  matrix into blocks
corresponds to the partition  $(k,1,k)$ of its rows and columns. Any
regular coadjoint orbit of the group $\UT$, where $n=2k+1$, contains
a unique element of the form
\begin{equation}\label{fff}f =
\begin{array}{lll}
\left(\begin{array}{ccc} 0&0&0\\
0&0&0\\
 \La&0&0\end{array}\right),& где& \La =
\left(\begin{array}{ccc}
0&\ldots&\la_k\\
\vdots&\ddots&\vdots\\
\la_1&\ldots&0\end{array}\right),
\end{array}\end{equation}
where $\la_1,\ldots,\la_{k}\in \Rb^*$ (see ~\cite[\S 3]{IP}).

\Theorem\Num \label{thodd}~(см. \cite{K1}). Let $n=2k+1$. The
character of the irreducible representation, associated with the
orbit of element  $f\in \ut^*$ of  (\ref{fff}), has the form

 \begin{equation}\label{regodd} \chi(g) = \de(C_{11}, C_{33}, C_{12}, C_{23})
 \frac{\chi_\La^*(C_{13})}{|\det \La|},\end{equation}

where $g$ as in  (\ref{ggg}), ~ $\chi_\La^*(C_{13})$ as in
(\ref{regeven}).

\Proof. Subalgebra
$ \left(\begin{array}{ccc} 0&0&*\\
0&0&*\\ 0&0&0\end{array}\right)$ is a polarization for  $f$. The
irreducible representation, associated with the orbit $\Omega(f)$,
is induced from the one-dimensional representation $e^{2\pi
i(\La,B_{13})}$ of the subgroup $\left\{ \left(\begin{array}{ccc} 0&0&B_{13}\\
0&0&B_{23}\\
0&0&0\end{array}\right)\right\}$. This representation is realized in
the space  $L^2(\Xx)$, where $\Xx$  consists of matrices
$$ X =
\left(\begin{array}{ccc} X_{11}&X_{12}&0\\
0&1&0\\
0&0&X_{33}
\end{array}\right),$$

by the formula

$$
T_gF(X) = e^{\dpi\left(\La,
(X_{11}C_{13}+X_{12}C_{23})(X_{33}C_{33})^{-1}\right)}F(X_{11}C_{11},
X_{11}C_{13} + X_{12}, X_{33}C_{33}).$$  Arguing as in subsection
2.1, we prove that $$ \chi(g) = \de(C_{11}, C_{33},C_{12}) \int
e^{\dpi \left(\La, (X_{11}C_{13} +
X_{12}C_{23})X_{33}\right)}dX_{11} dX_{12} dX_{33}. $$  Rewrite this
formula in the form

\begin{equation}\label{cccc}
\chi(g) = \de(C_{11}, C_{33}, C_{12})\cdot J\cdot \int e^{\dpi
\left(\La, X_{11}C_{13}X_{33}\right)}dX_{11} dX_{33},
\end{equation}
where $$ J = \int e^{\dpi \left(\La,
X_{12}C_{23}X_{33}\right)}dX_{12}.$$

  Let us show that  \begin{equation}\label{ii} J=\frac{1}{|\la_1\cdots
 \la_k|}\de(C_{23}).\end{equation}

Let
 $$
 X_{12} = \left(\begin{array}{l} x_{1,k+1}\\
 \vdots\\
 x_{k,k+1}\end{array}\right), \quad\quad C_{23} =
 (c_1,\ldots,c_k),\quad\quad X_{23} = \left(\begin{array}{cccc}
 1& y_{12}&\ldots&y_{1k}\\
0&1&\ldots&y_{2k}\\
\vdots&\vdots&\ddots&\vdots\\
0& 0&\ldots&1 \end{array}\right).
 $$
 Then $$\left(\La, X_{12}C_{23}X_{33}\right) = \sum_{s=1}^k \la_s M_s,$$
where
$$ M_s = x_{s,k+1}\left(c_1 y_{1,k-s+1} + \ldots + c_{k-s}
y_{k-s,k-s+1} + c_{k-s+1}\right).$$ In particular, $M_k =
x_{k,k+1}c_1$,~ $M_{k-1} = x_{k-1,k+1}\left(c_1y_{12}+c_2\right)$.
Denote
$$\left(\La, X_{12}C_{23}X_{33}\right)_t = \sum_{s=t}^k \la_s M_s,
\quad\quad J_t = \int e^{\dpi\left(\La, X_{12}C_{23}X_{33}\right)_t
}dx_{t,k+1}\cdots dx_{k,k+1}
$$
Note that
$$ \left(\La, X_{12}C_{23}X_{33}\right)_1= \left(\La,
X_{12}C_{23}X_{33}\right),\quad\quad J_1 = J,$$
$$ \left(\La, X_{12}C_{23}X_{33}\right)_t = \la_t M_t + \left(\La,
X_{12}C_{23}X_{33}\right)_{t+1}.$$ We use induction on $t$, moving
in decreasing order from $k$ to $1$, to prove
\begin{equation} \label{it}
J_t = \de(c_1,\ldots, c_{k-t+1})\frac{1}{|\la_t\cdots \la_k|}.
\end{equation} .

For $t=k$ we have
$$ J_k = \int e^{\dpi \la_kM_k}dx_{k,k+1} =  \int e^{\dpi \la_k
x_{k,k+1}c_1}dx_{k,k+1} = \frac{1}{|\la_k|}\de(c_1),
$$
this proves  (\ref{it}). Assume that  (\ref{it}) is proved for
$t+1$; let us prove it for  $t$:
$$J_t=  \int e^{\dpi\left(\La, X_{12}C_{23}X_{33}\right)_t
}dx_{t,k+1}\cdots dx_{k,k+1} =$$ $$  \int e^{\dpi
\la_tM_t}\left(\int e^{\dpi\left(\La,
X_{12}C_{23}X_{33}\right)_{t+1}}dx_{t+1,k+1}\cdots
dx_{k,k+1}dx_{t,k+1}\right) dx_{t,k+1}=$$
$$ \int e^{\dpi
\la_tM_t}\de(c_1,\ldots, c_{k-t})\frac{1}{|\la_{t+1}\cdots
\la_k|}dx_{t,k+1} = \de(c_1,\ldots, c_{k-t+1})\frac{1}{|\la_t\cdots
\la_k|}
$$
This proves  (\ref{it}). Substitute  (\ref{it}) into  (\ref{cccc}).
Calculation  of the integral  (\ref{cccc}) concludes similarly
subsection 2.1. $\Box$

\section{\textsc{Subregular characters}}

\subsection{Case of even  $n$.} Let $n=2(k+m+2)$. Partition rows and columns into blocks   $(m,1,1,k,k,1,1,m)$.
The general element  $g\in \UT$ can be written as a block matrix
\begin{equation}\label{geven}{\small g =
\left(\begin{array}{cccccccc} C_{11}&C_{12}&C_{13}&C_{14}&C_{15}&C_{16}&C_{17}&C_{18}\\
0&1&C_{23}&C_{24}&C_{25}&C_{26}&C_{27}&C_{28}\\
0&0&1&C_{34}&C_{35}&C_{36}&C_{37}&C_{38}\\
0&0&0&C_{44}&C_{45}&C_{46}&C_{47}&C_{48}\\
0&0&0&0&C_{55}&C_{56}&C_{57}&C_{58}\\
0&0&0&0&0&1&C_{67}&C_{68}\\
0&0&0&0&0&0&1&C_{78}\\
0&0&0&0&0&0&0&C_{88}\\
\end{array}\right).}
\end{equation}
It follows from  ~\cite[\S 3]{IP} that any subregular coadjoint
orbit of the group $\UT$, where $n=2(k+m+2)$, contains a unique
element
   of the form
   \begin{equation}\label{feven} {\small
f =
\left(\begin{array}{cccccccc} 0&0&0&0&0&0&0&0\\
0&0&0&0&0&0&0&0\\
0&0&0&0&0&0&0&0\\
0&0&0&0&0&0&0&0\\
0&0&0&\La_2&0&0&0&0\\
0&\ga_1&0&0&0&0&0&0\\
0&0&\ga_2&0&0&\ga_3&0&0\\
\La_1&0&0&0&0&0&0&0\\
\end{array}\right), }
\end{equation}
where $\La_1$ (resp. $\La_2$) is a quadratic matrix of the form
(\ref{ff}) of size  $m$ (resp. $k$). All  entries of the second
diagonal of   the matrix $\La_1$ and $\ga_1,\ga_2$  do not  equal to
zero. All  entries of the second diagonal of the matrix $\La_2$,
besides the last entry (see (\ref{ff})), do not equal to zero.

Denote by $E_{45}$ a system  of  $k$ matrix units that correspond to
the places of second diagonal of the block $C_{45}$ of the matrix
$g$ (see (\ref{geven})). Similarly, $E_{18}$  is a system of $k$
matrix units that correspond to the places of second diagonal of the
clock  $C_{18}$. The following lemma is proved by direct calculation.\\
 \Lemma \Num \label{ll}. 1)
Stabilizer
 $\gx^f$ in the Lie algebra  $\gx=\ut$ is a subspace spanned by
 \begin{equation}\label{ed}
 E_{45}, E_{18}, E_{26},
 E_{27}, E_{37}, \gamma_2 E_{23} + \gamma_1 E_{67};\end{equation}
 2) Stabilizer  $G^f$ in the Lie group  $G=\UT$ equals to   $E+\gx^f$.

 Construct the following table  of size  $8\times 8$ . By symbol "X"\, we mark all places
 $(a,b)$ such that $E_{a,b}$ rises in  (\ref{ed}).

\begin{center}
 {\small
\begin{tabular}{|p{0.2cm}|p{0.2cm}|p{0.2cm}|p{0.2cm}|p{0.2cm}|p{0.2cm}|p{0.2cm}|p{0.2cm}|}
\hline &$\bullet$ & &$\bullet$ &&$\bullet$&&X\rule{0pt}{0.5cm}\\
\hline  & &X & &&X&X&\rule{0pt}{0.5cm}\\
\hline &  & &$\bullet$ &&$\bullet$&X&$\bullet$\rule{0pt}{0.5cm}\\
\hline  &  & & &X&&&\rule{0pt}{0.5cm}\\
\hline & &  &  &&$\bullet$&&$\bullet$\rule{0pt}{0.5cm}\\
\hline & &  &  &&&X&\rule{0pt}{0.5cm}\\
\hline & &  &  &&&&$\bullet$\rule{0pt}{0.5cm}\\
\hline & &  &  &&&&\rule{0pt}{0.5cm}\\
\hline
\end{tabular}}\\
Table 1.
\end{center}

For any entry  $c$ of the matrix (\ref{geven}) that is lying in one
of the  blocks, marked in the Table 1 by the symbol  "$\bullet$"\,
or that is lying in the block  $C_{18}$, we define a rational
function $\ct$ as follows. Denote  $\Phi = \{(2,3),~(4,5),
~(6,7)\}$. Let $c$ be an entry of the block  $C_{ij}$  and
$\Phi(c)$ be a set of all pairs $(a,b)\in\Phi$ such that  $a > i$
and $b<j$. Let  $\De(c)$ be a minor of the matrix (\ref{geven}) with
rows and columns of $\Phi(c)$. Let  $\tilde{\De}(c)$ be a  minor of
the matrix (\ref{geven}), constructed by adding the row (resp.
column) of the entry $c$ to the system of rows (resp. columns) of
the minor $\De(c)$. Denote
  \begin{equation}\label{ct}
  \ct = \pm \De(c)^{-1} \tilde{\De}(c),\end{equation}
  where the sign  $\pm$ coincides with the sign of entry  $c$ in the minor
  $\tilde{\De}(c)$.
  For example,  if $m=k=1$ (i.e. every block  $C_{ij}$ is an entry
  $c_{ij}$), then\\
  \\$ \ct_{12} = c_{12}$,~ $ \ct_{34} = c_{34}$,~ $ \ct_{56} =
  c_{56}$,~   $ \ct_{78} = c_{78}$,~ $\ct_{14} = - c_{23}^{-1}
  \left|\begin{array}{cc} c_{13}& c_{14}\\
  c_{23}& c_{23}\end{array}\right|$, \\
 $\ct_{36} = - c_{45}^{-1}
  \left|\begin{array}{cc} c_{35}& c_{36}\\
  c_{45}& c_{46}\end{array}\right|$, ~ $\ct_{58} = - c_{67}^{-1}
  \left|\begin{array}{cc} c_{57}& c_{58}\\
  c_{67}& c_{68}\end{array}\right|$,~
  $\ct_{16} =  c_{23}^{-1} c_{45}^{-1}
  \left|\begin{array}{ccc} c_{13}& c_{15}& c_{16}\\
  c_{23}& c_{25}& c_{26}\\
  0& c_{45}& c_{46}\end{array}\right|$, ~
$\ct_{38} =  c_{45}^{-1} c_{67}^{-1}
  \left|\begin{array}{ccc} c_{35}& c_{37}& c_{38}\\
  c_{45}& c_{47}& c_{48}\\
  0& c_{67}& c_{68}\end{array}\right|$,
  ~~$
  \ct_{18} = -c_{23}^{-1} c_{45}^{-1} c_{67}^{-1}
  \left|\begin{array}{cccc} c_{13}& c_{15}& c_{17}& c_{18}\\
  c_{23}& c_{25}& c_{27}& c_{28}\\
  0&c_{45}& c_{47}& c_{48}\\
  0&0& c_{67}& c_{68}\end{array}\right|.$\\
Denote: \\
1)~ $S_0$ is the set of entries lying upper the diagonal in the
blocks  $
C_{11},C_{44}, C_{55}, C_{88}$;\\
 2)~ $S_1$ is the set of all rational functions  $\ct$, where $c$
 runs through all entries in the blocks that is marked in the Table 1 by the symbol  "$\bullet$"\,;\\
 3) $S_2 = \{ C_{23}C_{35} + C_{24}C_{45}, ~ C_{45}C_{57} + C_{46}C_{67}, ~
 \gamma_1C_{23} - \gamma_2
 C_{67}\}$,
  \\
 4)~ $S = \{S_0, S_1, S_2\}$,\\
5)~   $d(g) =  \det C_{23}\cdot\det C_{45}\cdot \det C_{67}$. Note
that according to partition  into blocks  (\ref{geven}), the blocks
 $C_{23}$,~ $C_{67}$ have size 1 and
 $\det C_{23} = C_{23}$, ~ $\det C_{67} = C_{67}$.

 The algebra   $\Rb[G]$ of regular functions on the group
 $G=\UT$ admits localization  $\Rb'[G]$ by the denominator system generated by $d(g)$.
Note that  $S\subset \Rb'[G]$. Let  $I'$ be an ideal in  $\Rb'[G]$
generated by  $S$. The annihilator of the ideal  $I'$ is the subset
in $G' = \{g\in G: ~ d(g)\ne 0\}.$\\
\Prop\Num \label{ii}. Let   $f$ be of the form  (\ref{feven}). Then
the closure in $G$ of the  annihilator of the ideal  $I'$ coincides
with the closure of the set  $\Ad_G(G^f)$.\\
\Proof. The Lemma \ref{ll} implies that a function of  $S$
annihilate  $G^f$. One can prove directly that  $S$ annihilate
$\Ad_G(G^f)$. Then  $\Ann I'\supset \Ad_G(G^f)\cap G'$. The ideal
$I'$ is prime and its dimension coincides with the dimension of the
set   $\Ad_G(G^f)$. $\Box$
\\
 \Theorem\Num\label{subeven}. Let  $n=2(k+m+2)$.
The character of the irreducible representation, associated with the
orbit of element  $f\in \ut^*$ of  (\ref{feven}), has the form
\begin{equation}\label{cheven}
\chi(g) = \de(S)\cdot
\chi_{\La_1}^*(\widetilde{C}_{18})\cdot\chi_{\La_2}^*(C_{45})\cdot\chi_0^*(g)
\end{equation}
where $\chi_{\La_2}^*(C_{45})$  as in  (\ref{regeven}),~ the matrix
$\widetilde{C}_{18}$ is filled by  entries  $\ct$, where $c$ is an
corresponding entry of the block  $C_{18}$ (see (\ref{ct}),~
$\chi_{\La_1}^*(\widetilde{C}_{18})$ as in  (\ref{regeven}), and
$$\chi_0^*(g) = \frac{1}{|\ga_1\ga_2|^k \cdot d(g)^m}
e^{\dpi \left(\frac{\gamma_1}{C_{67}}Q_0 + \ga_3 C_{67}\right)},
$$
$$
Q_0= c_{23} c_{37} + c_{24} c_{47} + c_{25} c_{57} + c_{26}
 c_{67}.$$
Remark,  since the element  $\gamma_1C_{23} - \gamma_2 C_{67}$
belongs to
 the set  $S$, one can substitute $\frac{\gamma_1}{C_{67}}$ by
$\frac{\gamma_2}{C_{23}}$ in the formula for $\chi_0^*(g)$.\\
{\bf Proof scheme}.
\\
{\bf Item 1}. Consider the special case $m=0$. Partition of rows and
columns has the form  $(1,1,k,k,1,1)$. The general element  $g\in
G=\UT$ and  $f_0\in \gx^*$ are written as block matrices
$${\small
g =
\left(\begin{array}{cccccc} 1&C_{12}&C_{13}&C_{14}&C_{15}&C_{16}\\
0&1&C_{23}&C_{24}&C_{25}&C_{26}\\
0&0&C_{33}&C_{34}&C_{35}&C_{36}\\
0&0&0&C_{44}&C_{45}&C_{46}\\
0&0&0&0&1&C_{56}\\
0&0&0&0&0&1\\
\end{array}\right), \quad\quad\quad f_0 =
\left(\begin{array}{cccccc} 0&0&0&0&0&0\\
0&0&0&0&0&0\\
0&0&0&0&0&0\\
0&0&\La&0&0&0\\
\ga_1&0&0&0&0&0\\
0&\ga_2&0&0&\ga_3&0\\
\end{array}\right)}
$$

Consider the subgroup $B$ of all matrices
$${\small b = \left(\begin{array}{cccccc} 1&b_{12}&b_{13}&b_{14}&b_{15}&b_{16}\\
0&1&0&0&0&b_{26}\\
0&0&b_{33}&b_{34}&0&b_{36}\\
0&0&0&b_{44}&0&b_{46}\\
0&0&0&0&1&b_{56}\\
0&0&0&0&0&1\\
\end{array}\right)}
$$
Let  $\pi_\La(b_1)$ be the irreducible representation from
subsection 2.1 of the unitriangular group  $B_1$ of all  matrices
$b_1 =
\left(\begin{array}{cc} b_{33}& b_{34}\\
0& b_{44}\end{array}\right)$

 The irreducible representation  $T^{f_0}_g$, that corresponds to $f_0\in\gx^*$,
 is induced from the representation  $$e^{\dpi
\left(\ga_1b_{15}+\ga_2 b_{26}+\ga_3b_{56}\right)}\pi_\La(b_1)$$ of
the subgroup  $B$. Direct calculations in spirit of subsection  2.1
lead to proof of formula (\ref{cheven}) for special case  $m=0$.

{\bf Item  2}. Consider the special case $m=1$. Partition of rows
and columns has the form  $(1,1,1,k,k,1,1,1)$. The general element
$g\in G=\UT$ and $f\in \gx^*$ are written as block matrices

$${\small g =
\left(\begin{array}{cccccccc} 1&C_{12}&C_{13}&C_{14}&C_{15}&C_{16}&C_{17}&C_{18}\\
0&1&C_{23}&C_{24}&C_{25}&C_{26}&C_{27}&C_{28}\\
0&0&1&C_{34}&C_{35}&C_{36}&C_{37}&C_{38}\\
0&0&0&C_{44}&C_{45}&C_{46}&C_{47}&C_{48}\\
0&0&0&0&C_{55}&C_{56}&C_{57}&C_{58}\\
0&0&0&0&0&1&C_{67}&C_{68}\\
0&0&0&0&0&0&1&C_{78}\\
0&0&0&0&0&0&0&1\\
\end{array}\right), \quad\quad f =
\left(\begin{array}{cccccccc} 0&0&0&0&0&0&0&0\\
0&0&0&0&0&0&0&0\\
0&0&0&0&0&0&0&0\\
0&0&0&0&0&0&0&0\\
0&0&0&\La_2&0&0&0&0\\
0&\ga_1&0&0&0&0&0&0\\
0&0&\ga_2&0&0&\ga_3&0&0\\
\la_1&0&0&0&0&0&0&0\\
\end{array}\right), }
$$ where $\la_1\in \Rb^*$.

Consider the subgroup  $B_*$ of all matrices
$${\small b_* = \left(\begin{array}{cccccccc} 1&0&0&0&b_{15}&b_{16}&b_{17}&b_{18}\\
0&1&b_{23}&b_{24}&b_{25}&b_{26}&b_{27}&b_{28}\\
0&0&1&b_{34}&b_{35}&b_{36}&b_{37}&b_{38}\\
0&0&0&b_{44}&b_{45}&b_{46}&b_{47}&b_{48}\\
0&0&0&0&b_{55}&b_{56}&b_{57}&0\\
0&0&0&0&0&1&b_{67}&0\\
0&0&0&0&0&0&1&0\\
0&0&0&0&0&0&0&1\\
\end{array}\right)}
$$

Let  $b_*'$ be a matrix obtained from  $b_*$ by deleting
 the first row and the last column, $B_*'$  be a group of all such matrices.
 Let  $T^{f_0}(b_*')$ be an irreducible representation of the group $B_*'$ as in Item 1.

 The irreducible representation $T^f_g$, that corresponds to
$f\in\gx^*$, is induced from the representation
$$e^{\dpi \left(\la_1b_{18}\right)}T^{f_0}(b_*)$$ of the subgroup $B$.
Direct calculations in spirit of subsection  2.1 lead to proof of
formula (\ref{cheven}) for special case  $m=1$.
\\
{\bf Item 3}. Case of general  $m$. Consider the normal subgroup
$B_m$ of all  $g\in G$, such that  $C_{11}=C_{88}=E$. The character
 $\chi_m$  of irreducible representation of the subgroup  $B_m$,
that corresponds to the restriction of  $f$ on $\mathrm{Lie}(B_m)$,
is calculated similarly to Item  2.
 The character  $\chi$ is induced from the character  $\chi_m$ of subgroup $B_m$. Calculation of character is similar
 to subsection  2.1. $\Box$

\subsection{Case of odd  $n$.} Let  $n=2(k+m+2)+1$.
Partition rows and columns in blocks   $(m,1,1,k,1,k,1,1,m)$. The
general element  $g\in \UT$ can be written as a block matrix
\begin{equation}\label{godd}{\small g =
\left(\begin{array}{ccccccccc} C_{11}&C_{12}&C_{13}&C_{14}&C_{15}&C_{16}&C_{17}&C_{18}&C_{19}\\
0&1&C_{23}&C_{24}&C_{25}&C_{26}&C_{27}&C_{28}&C_{29}\\
0&0&1&C_{34}&C_{35}&C_{36}&C_{37}&C_{38}&C_{39}\\
0&0&0&C_{44}&C_{45}&C_{46}&C_{47}&C_{48}&C_{49}\\
0&0&0&0&1&C_{56}&C_{57}&C_{58}&C_{59}\\
0&0&0&0&0&C_{66}&C_{67}&C_{68}&C_{69}\\
0&0&0&0&0&0&1&C_{78}&C_{79}\\
0&0&0&0&0&0&0&1&C_{89}\\
0&0&0&0&0&0&0&0&C_{99}
\end{array}\right).}
\end{equation}
   Any subregular coadjoint orbit of the group $\UT$, where $n=2(k+m)+5$, contains a unique element
   of the form
   \begin{equation}\label{fodd} {\small
f =
\left(\begin{array}{ccccccccc} 0&0&0&0&0&0&0&0&0\\
0&0&0&0&0&0&0&0&0\\
0&0&0&0&0&0&0&0&0\\
0&0&0&0&0&0&0&0&0\\
0&0&0&0&0&0&0&0&0\\
0&0&0&\La_2&0&0&0&0&0\\
0&\ga_1&0&0&0&0&0&0&0\\
0&0&\ga_2&0&0&0&\ga_3&0&0\\
\La_1&0&0&0&0&0&0&0&0\\
\end{array}\right), }
\end{equation}
where $\La_1$ (resp. $\La_2$)  is a matrix of form  (\ref{ff}) of
size $m$ (resp. $k$). All  entries of second diagonals of the
matrices $\La_1$,~$\La_2$ and $\ga_1,\ga_2$ do not equal to zero
~\cite[\S 3]{IP}).

The systems of matrix units  $E_{46}$,~ $E_{19}$ is defined by
blocks
$C_{46}$,~ $C_{19}$ similarly as in subsection 3.1.\\
 \Lemma \Num \label{lll}.
 1)
Stabilizer  $\gx^f$ in the Lie algebra $\gx=\ut$ is a subspace
spanned by
 \begin{equation}\label{edd}
 E_{46}, E_{19}, E_{27},
 E_{28}, E_{38}, \gamma_2 E_{23} + \gamma_1 E_{78};\end{equation}
 2) Stabilizer  $G^f$ in the Lie group  $G=\UT$ is equal to   $E+\gx^f$.

Construct the following table  of size  $9\times 9$ . By symbol
"X"\, we mark all places
 $(a,b)$ such that $E_{a,b}$ rises in  (\ref{edd}).

\begin{center}
 {\small
\begin{tabular}{|p{0.2cm}|p{0.2cm}|p{0.2cm}|p{0.2cm}|p{0.2cm}|p{0.2cm}|p{0.2cm}|p{0.2cm}|p{0.2cm}|}
\hline &$\bullet$ & &$\bullet$ &$\bullet$&&$\bullet$&&X\rule{0pt}{0.5cm}\\
\hline  & &X & &&&X&X&\rule{0pt}{0.5cm}\\
\hline &  & &$\bullet$ &$\bullet$&&$\bullet$&X&$\bullet$\rule{0pt}{0.5cm}\\
\hline  &  & & &$\bullet$&X&&&\rule{0pt}{0.5cm}\\
\hline & &  &  &&$\bullet$&$\bullet$&&$\bullet$\rule{0pt}{0.5cm}\\
\hline & &  &  &&&$\bullet$&&$\bullet$\rule{0pt}{0.5cm}\\
\hline & &  & & &&&X&\rule{0pt}{0.5cm}\\
\hline & &  & & &&&&$\bullet$\rule{0pt}{0.5cm}\\
\hline & &  & & &&&&\rule{0pt}{0.5cm}\\
\hline
\end{tabular}}\\
Table  2.
\end{center}

Denote $\Phi' = \{(2,3),~(4,6), ~(7,8)\}$. For any entry  $c$ of the
matrix (\ref{geven}) that is lying in one of the  blocks, marked in
the Table 2 by the symbol  "$\bullet$"\,, or that is lying in the
block  $C_{18}$, we define a rational function $\ct$ as in
(\ref{ct}) changing $\Phi$ by $\Phi'$.
 Denote: \\
1)~ $S_0$ is the set of entries lying upper the diagonal in the
blocks $
C_{11},C_{44}, C_{66}, C_{99}$;\\
 2)~ $S_1$ is the set of all rational functions  $\ct$, where $c$
 runs through all entries in the blocks that is marked in the Table 2 by the symbol  "$\bullet$"\,;\\
 3) $S_2 = \{ C_{23}C_{36} + C_{24}C_{46}, ~ C_{46}C_{68} + C_{47}C_{78}, ~ \gamma_1C_{23} - \gamma_2
 C_{78}\}$,
  \\
 4)~ $S = \{S_0, S_1, S_2\}$,\\
5)~  $d_1(g) =  \det C_{23}\cdot\det C_{46}\cdot \det C_{78}$.

 The
algebra   $\Rb[G]$ of regular function on the group $G=\UT$ admits
localization $\Rb'[G]$ by the denominator subset generated by
$d_1(g)$. As above denote by  $I'$ the ideal in $\Rb'[G]$ generated
by $S$; note that   $\Ann(I')\subset \{g\in
G: ~ d(g)\ne 0\}.$\\
\Prop\Num \label{iii}. Let   $f$ has the form (\ref{fodd}).
  Then
the closure in $G$ of the  annihilator of the ideal  $I'$ coincides
with the closure of the set  $\Ad_G(G^f)$.\\
\Proof. Similarly to Proposition \ref{ii}.
\\
 \Theorem\Num \label{subodd}. Let $n=2(k+m)+5$.
The character of the irreducible representation, associated with the
orbit of element  $f\in \ut^*$ of  (\ref{fodd}), has the form
\begin{equation}\label{codd}
\chi(g) = \de(S)\cdot
\chi_{\La_1}^*(\widetilde{C}_{19})\cdot\frac{\chi_{\La_2}^*(C_{46})}{|\det
\La_2|}\cdot\chi_1^*(g).
\end{equation}
Here $\chi_{\La_2}^*(C_{46})$  as in  (\ref{regeven}), the matrix
$\widetilde{C}_{19}$  filled by the entries  $\ct$, where $c$ is the
corresponding entry of the block  $C_{19}$ (see (\ref{ct}),~
$\chi_{\La_1}^*(\widetilde{C}_{19})$ as in  (\ref{regeven}), and
$$\chi_1^*(g) = \frac{1}{|\ga_1\ga_2|^k \cdot|\ga_1c_{23}|\cdot d_1(g)^m}
e^{\dpi \left(\frac{\gamma_1}{C_{67}}Q_1 + \ga_3 C_{67}\right)},
$$
$$
Q_1 = c_{23} c_{38} + c_{24} c_{48} + c_{25} c_{58} + c_{26} c_{68}
+ c_{27}c_{78}.$$
\Proof ~ is similar to Theorem  \ref{subeven}. \\
\Theorem\Num \label{thhy}. Conjecture  \ref{hy} is valid for regular
and subregular characters of the group  $\UT$.
\\
\Proof.~ Propositions \ref{ii}, \ref{iii} and Theorems \ref{theven},
\ref{thodd}, \ref{subeven}, \ref{subodd} imply the proof.

\end{document}